\documentclass{article}
\usepackage[utf8]{inputenc}
\usepackage{amsmath,amsfonts,amssymb}
\usepackage{amsthm}
\usepackage[noadjust]{cite}
\newcommand{\e}{\mathsf{E}}
\newcommand{\p}{\mathsf{P}}
\newcommand{\di}{\mathsf{Var}}

\newtheorem*{cond1}{Condition}
\newtheorem{thm1}{Theorem}

\newtheorem{lm1}{Lemma}
\newtheorem{rem1}{Remark}

\title{Strong Gaussian approximation for cumulative processes with heavy tails}
\author{Elena Bashtova\footnote{Lomonosov Moscow State University; supported by RFBR grant 20-01-00487; e-mail: elena.bashtova@math.msu.ru}, Alexey Shashkin }
\date{}

\begin{document}

\maketitle

\section{Introduction}

This paper is a continuation of work \cite{ASAExp} devoted to establishment of the convergence rate in the strong invariance principle for cumulative processes. Cumulative processes were introduced by Smith \cite{Smith1955} and attract much attention in queueing systems theory, see discussion in \cite{ASAExp}. Here we consider the case when the regeneration periods and increments over them have the moments of order $p>2.$ It is well-known (see, e.g., the survey \cite{Zait_survey}) that under power-type conditions one can be interested in two types of approximation by etalon Wiener process: the rate of convergence in the Strassen's invariance principle and the estimation of probability that random process deviates from the approximating one. 
That is, given a stationary centered sequence $\{X_i,i\geq 1\}$ and setting $S_n=X_1+\ldots+X_n,$ $n\geq 1,$ we can try to construct the approximating scaled Wiener process $\{\sigma W_t,t\geq 0\}$ ($\sigma^2$ being the asymptotic variance of $S_n$) such that either
\begin{equation}\label{intr1}\left|S_n -\sigma W_n\right| = o(f(n)) \mbox{ a.s.} \end{equation}
for some function $f:\mathbb{N}\to\mathbb{R}_+$ growing slowly enough, 
or \begin{equation}\label{intr2}\p\left(\sup_{k\leq n}|S_k-\sigma W_k|\geq x\right)\leq a_0nx^{-p} \end{equation}
for a non-random $a_0>0$ and $x$ lying in some domain depending on the value of  $t.$ Note that when the $X_i$ are i.i.d. satisfying the Cram\'{e}r condition, one can obtain the exponential inequality of the type \eqref{intr2} which implies \eqref{intr1} with $f(n)\sim \log n.$ However under the conditions of the type $\e H(X_1)<\infty$ with even and regularly varying $H,$ these two properties typically need to be established separately. We confine ourselves to the case $H(x)=|x|^p$ with $p>2,$ since in that case one can hope to get the non-improvable rate \eqref{intr1} with $f(n)\sim n^{1/p}.$

 We study a $\mathbb{R}^d-$valued random process $S=\{S(t),t\geq 0\} = \{(S_1(t),\ldots,S_d(t)),t\geq 0\}$ which is assumed to be separable. For a vector $z\in\mathbb{R}^d$, by $|z|$ we denote the maximal norm. The principal condition we will require is: 
 
 \begin{cond1}[{\textbf{A}}]
 There exists an a.s. increasing random sequence $\{T_k, k\geq 0\}$ such that $T_0=0$ and random elements
 $$ \Bigl\{\bigl(T_j-T_{j-1}, S(T_{j-1}+t) - S(T_{j-1}), t\in(0,T_j-T_{j-1}]  \bigr) ,j\geq 1\Bigr\}$$
 are i.i.d.
 \end{cond1}

 We need the following notation. 
\begin{itemize} 
\item $\tau_k = T_k - T_{k-1},\;k\geq 1;$ \quad $\xi_k =(\xi_{k1},\ldots,\xi_{kd})= S(T_k) - S(T_{k-1}),\; k\geq 1;$
\item  $m(t) = \max\limits\{k\in\mathbb{Z}_+: T_k \leq t\},\; t\geq 0$, i.e. $m=\{m(t),t\geq 0\}$ is the renewal process built by sequence $\{\tau_k,k\geq 1\}$;
\item $\eta_k = \max\limits_{0< t \leq T_k-T_{k-1}}|S(T_{k-1}+t) - S(T_{k-1})|,\;k\geq 1;$ 
\item $\mu=\e \tau_1$, \quad $\varkappa = \mu^{-1}\e \xi_1  ;$ \quad  $\sigma^2 =   \mu^{-1}(\di(\xi_1) - 2 cov(\xi_1,\tau_1)\varkappa^T + \varkappa \varkappa^T \di(\tau_1)),$
\item  \begin{equation}\label{greeks}
\beta = \bigl(\di (\tau_1)\bigr)^{-1} {cov(\xi_1, \tau_1)},\; v^2 = \di(\xi_1 - \beta \tau_1),\;
\gamma = \frac{\di(\tau_1)}{\mu}, \;\lambda = \frac{\mu^2 }{\di(\tau_1)},\;\alpha = \beta - \varkappa.
\end{equation}
 \end{itemize}
 where $cov(\xi_1,\tau_1) := (cov(\xi_{11},\tau_1),\ldots, cov(\xi_{1d},\tau_1) ).$ We write $\sigma$ for the usual square root of the matrix $\sigma^2.$

 Next we invite the moment  condition.
 \begin{cond1}[\textbf{Bp}]
There is $p>2$ such that $\e  \tau_1^p <\infty$ and $\e  \eta_1^p <\infty$.
 \end{cond1}

\section{Main results}

\begin{thm1}\label{TheoremMain2} 
Suppose that $\{S(t),t\geq 0\}$ satisfies conditions {\bf(A)} and {\bf{(Bp)}}.  Then one can redefine it on a probability space $(\Omega,\mathcal{F},\mathsf{P})$ supporting also a standard $d$-dimensional Wiener process $\{W_t,t\geq 0\}$ such that the relation
\begin{equation}\label{1t}
 \sup_{u\leq t} \bigl|S(u) - \varkappa u - \sigma W_u\bigr| =  o(t^{1/p})
\end{equation}
holds as $t\to\infty$ with probability one.
\end{thm1}

\begin{thm1}\label{TheoremMain} 
Suppose that $\{S(t),t\geq 0\}$ satisfies conditions {\bf(A)} and {\bf{(Bp)}}.  Then one can redefine it on a probability space $(\Omega,\mathcal{F},\mathsf{P})$ supporting also a standard $d$-dimensional Wiener process $\{W_t,t\geq 0\}$ such that for some positive constant $a$   and all $t\geq 1, x>0$ one has
\begin{equation}\label{1}
\p\left(\sup_{u\leq t}\bigl|S(u) - \varkappa u - \sigma W_u\bigr| \geq x \right) \leq a t x^{-p} . 
\end{equation}
\end{thm1}

\begin{proof} We will start proving both theorems simultaneously, discussing the differences in the proofs where it is appropriate. While working on both types of strong approximation, we will use the same notation for a Wiener process approximating a given sequence in both senses, even though it cannot be claimed that assertions of both types can be satisfied using one process only. Also, clearly we can consider $t\geq e.$ Without loss of generality we can make the following assumptions:

\begin{itemize}
    \item 
  $\di (\tau_1)>0$ as otherwise the statements would be easily derived from the Koml\'{o}s-Major-Tusn\'{a}dy results for i.i.d. summands \cite[Theorems 2 and 4]{KMT2};
  \item 
for proving Theorem \ref{TheoremMain} assume that $x\geq ct^{1/p},$ with $c>0$ a  fixed constant to be chosen later, as for other $x$ the estimate is trivial. 
  \end{itemize}
We will use  positive constants $a_i,\,i\geq 1,$ when we want to write an auxiliary estimate of the type \eqref{1}.  


 
Note that the choice of $\beta$ implies that $cov(\xi_{ki} - \beta_i \tau_k,\tau_k)=0, \;i=1,\ldots,d. $ Applying Einmahl's \cite{Einmahl,Berger} and Zaitsev's \cite{Zait_survey} theorems to a $(d+1)$-dimensional sequence 
$\{(\xi_k - \beta \tau_k +\alpha\mu,\tau_k-\mu), k\geq 1\}$ we can infer (redefining the initial process on a larger probability space if necessary) that there exist two independent Wiener processes $\{B_t,t\geq 0\},\{\widetilde{B}_t,t\geq 0\} $, the former being $d-$dimensional, such that for any pair $(t,x)\in [1,\infty)\times (0,\infty)$

\begin{equation}\label{2}
\sup_{k\leq t}\bigl|S(T_k) + \alpha k  \mu - \beta T_k   - vB_k \bigr|=o(t^{1/p}),\qquad\sup_{k\leq t}\bigl|T_k - k\mu  - \sqrt{\di(\tau_1)}\widetilde{B}_k \bigr|=o(t^{1/p})\mbox{ a.s.}, 
\end{equation}
\begin{equation}\label{3}
\p\left(\sup_{k\leq t}\bigl|S(T_k) + \alpha k  \mu - \beta T_k   - vB_k \bigr| \geq x \right) \leq  \frac{a_1t} {x^{p}},
\quad\p\left(\sup_{k\leq t}\bigl|T_k - k\mu  - \sqrt{\di(\tau_1)}\widetilde{B}_k \bigr| \geq x \right) \leq  \frac{a_1t}{ x^{p}} , 
\end{equation}
with some $a_1>0.$

\begin{rem1}\label{Rem1} Note that given a random process $\{Z_t,t\geq 0\}$ and a bound
$$\p\left(\sup_{u\leq t}|Z_u|\geq C\log t+z\right)\leq Ae^{-Bz} $$
holding for some $A,B,C>0$ and all $t\geq 1,z\geq 0,$ by adjusting the factor $c$ (namely, letting $c\geq Cp$) we can always infer that an estimate like in \eqref{1} is valid for $x\geq c t^{1/p}.$ Also, same bound implies that $\sup_{u\leq t}|Z_u|$ is $O(\log t)$ as $t\to\infty,$ with probability one.
\end{rem1}
Consequently, applying formula (5) from \cite{ASAExp} (following the idea of \cite{MerRio}), we can construct Poisson process $\{N(t),t\geq 0\}$ with parameter $\lambda$, measurable with respect to the $\sigma-$algebra generated by $\widetilde{B}$ and such that the following relations hold:

\begin{equation}\label{5}
\sup_{u\leq t}\bigl|\gamma N_u - u\mu- \sqrt{\di(\tau_1)}\widetilde{B}_u \bigr|=o(t^{1/p}), \qquad \sup_{u\leq t}\bigl|\gamma N_u - T_{[u]} \bigr| =o(t^{1/p})\mbox{ a.s.},  
\end{equation}
\begin{equation}\label{4a}
\p\left(\sup_{u\leq t}\bigl|\gamma N_u - u\mu- \sqrt{\di(\tau_1)}\widetilde{B}_u \bigr| \geq x \right) \leq a_2 t x^{-p} ,\quad \p\left(\sup_{u\leq t}\bigl|\gamma N_u - T_{[u]} \bigr| \geq x \right) \leq a_2 t x^{-p}  
\end{equation}
with some $a_2>0.$

Further, we denote $y(u) = N^{-1}(u/\gamma)$. Then, process $m$ is the inverse of process $T$ and $y$ is the inverse of   $\gamma N$. Therefore in view of the paper of  Cs\"{o}rg\H{o}, Horv\'{a}th and Steinebach \cite{CHS} and particular Theorems 3.2 and 4.1 there, and of \eqref{2}--\eqref{4a}, one can construct (enlarging the probability space if necessary) a standard Wiener process $\widetilde{W}=\{\widetilde{W}_t,t\geq 0\}$ such that:

\begin{equation}\label{6}
 \sup_{u\leq t}\left|m(u) - \frac{u}{\mu} - \frac{1}{\lambda\sqrt{\gamma}}\widetilde{W}_u\right|=o(t^{1/p}),\qquad  
 \sup_{u\leq t}\left| y(u)-  \frac{u}{\mu} - \frac{1}{\lambda\sqrt{\gamma}}\widetilde{W}_u \right|=o(t^{1/p})\mbox{ a.s.},
\end{equation}
\begin{equation}\label{6a}
 \sup_{u\leq t}  |m(u) - {y(u)}|=o(t^{1/p})\mbox{ a.s.},
\end{equation}
\begin{equation}
\label{deheuvels}
\p\left(\sup_{u\leq t}\left|m(u) - \frac{u}{\mu} - \frac{1}{\lambda\sqrt{\gamma}}\widetilde{W}_u\right|\geq   x\right)  \leq a_3 t x^{-p},\quad \p\left(\sup_{u\leq t}\left| y(u)-  \frac{u}{\mu} - \frac{1}{\lambda\sqrt{\gamma}}\widetilde{W}_u \right| \geq  x \right) \leq a_3 t x^{-p} ,
\end{equation}
and
\begin{equation}\label{m-y}
\p\left(\sup_{u\leq t}  |m(u) - {y(u)}| \geq   
  x \right) \leq a_3 t x^{-p},
 \end{equation}
 with some $a_3>0$ and $x\in[ct^{1/p}, t/\log t].$ Here, depending on which of the Theorems \ref{TheoremMain2} and \ref{TheoremMain} is being proved, either \eqref{6}--\eqref{6a} or \eqref{deheuvels}--\eqref{m-y} are simultaneously true.  Moreover, process $\widetilde{W}$ is measurable with respect to a $\sigma$-algebra generated by $\widetilde{B}$ plus some random element $V$ independent from all ones considered above.

 By Lemmas 2 and 3 in \cite{ASAExp} and Remark \ref{Rem1}, one can construct a standard $d$-dimensional Wiener process $W^*=\{W^*_t,t\geq 0\}$ such that  
\begin{equation}\label{7a}\sup_{u\leq t}\left|  B_u - \frac{W^*_{N(u)}}{\sqrt{\lambda}}   \right|  =O(\log t)\mbox{ a.s.,}  \end{equation}
\begin{equation}\label{7}\p\left(\sup_{u\leq t}\left|  B_u - \frac{W^*_{N(u)}}{\sqrt{\lambda}}   \right| \geq   x \right) \leq a_4 t x^{-p}  \end{equation}
with some $a_4>0$ and all $t\geq e,\,x \geq ct^{1/p} ;$ this process is determined by the process $B$ plus some random element $V^*$ which is independent from all other random elements considered above. In particular, $W^*$  is independent from $\widetilde{B}$ and $N$.

 By construction,  standard Wiener processes $\{\widetilde{W}_t,t\geq 0\}, \{W^*_t,t\geq 0\}$, the latter being $d$-dimensional,  are independent.

Let $W^{\circ}=\{W^{\circ}_t,t\geq 0\}$ be a standard $d$-dimensional Wiener process independent from $(W^*,\widetilde{W}).$ Define a Gaussian   process
\begin{equation}\label{4}
W_t = \sigma^{+}\left( \lambda^{-1/2}{v} W^*_{t/\gamma}  - \lambda^{-1}\gamma^{-1/2} \mu \alpha  \widetilde{W}_{t} \right) + (I_d - \sigma^{+} \sigma)W^{\circ}_t, \end{equation}
where $v$ is the square root of the covariance matrix $v^2,$ and $\sigma^+$ stands for the Moore-Penrose pseudo-inverse matrix to $\sigma.$
By Lemma 4 in \cite{ASAExp} $W$  is a standard $d$-dimensional Wiener process.

Denoting $y = y(u) = N^{-1}(u/\gamma) $ and using relations $\sigma\sigma^+ v = v ,$  $\sigma\sigma^+ \alpha=\alpha $  we write, for $u\geq 0,$
\begin{multline}\label{eight_phis}
S(u) - \varkappa u -\sigma W_u = S(u)   - \varkappa u - \sigma\sigma^+\frac{v}{\sqrt{\lambda}}W^*_{u/\gamma} + \sigma\sigma^+\alpha \frac{\mu}{\lambda\sqrt{\gamma}}\widetilde{W}_{u} =\\=S(u)   - \varkappa u - \frac{v}{\sqrt{\lambda}}W^*_{u/\gamma}  -\alpha \mu \left(N^{-1}(u/\gamma) - \frac{u}{\lambda \gamma} - \frac{\widetilde{W}_{u }}{\lambda\sqrt{\gamma}}\right) + \alpha\mu\left(N^{-1}(u/\gamma) - \frac{u}{\lambda \gamma} \right) =
\\= (S(u) - S(T_{m(u)}) )+   (S(T_{m(u)}) - S(T_{[y]})) +\left(S(T_{[y]}) -\beta T_{[y]}   + \alpha \mu y - v B_y\right)+ \\ +\beta  (T_{[y]} - \gamma N(y)) - \alpha \mu \left(y - \frac{u}{\lambda \gamma} - \frac{\widetilde{W}_{u }}{\lambda\sqrt{\gamma}}\right)  +v\left(B_y - \lambda^{-1/2}W^*_{N(y)}\right)+\\+\lambda^{-1/2}v\left(  W^*_{N(y)}-  W^*_{u/\gamma}\right)+\beta(\gamma[u/\gamma]+\gamma-u)  =:\sum_{q=1}^8 \Phi_q(u).  \end{multline}
 
Now we will proof Theorem \ref{TheoremMain2}. Note that by strong law of large numbers $m(t)\sim t/\mu$ and $N^{-1}(t/\gamma)\sim t/\mu$ a.s. when $t\to\infty.$  We need to show that $\sup_{u\leq t}|\Phi_q(u)| = o(t^{1/p})$ a.s. for each $q=1,\ldots,8.$ For $q=1$ this follows from the estimate 
$$\limsup_{t\to\infty}t^{-1/p}\sup_{u\leq t}|S(u) - S(T_{m(u)})|\leq \limsup_{t\to\infty}t^{-1/p}\max_{k\leq 2t/\mu }
\eta_k = 0
$$
by the condition {\bf(Bp)} and standard bound for maxima, see e.g. Embrechts,  Kl\"{u}ppelberg and Mikosch \cite[Theorem 3.5.1]{EMK}.
For $q=2,$    using \eqref{2} and \eqref{6a}, for any fixed $\varepsilon>0$ we get
$$\limsup_{t\to\infty}t^{-1/p}|\Phi_2(t)|\leq \limsup_{t\to\infty}t^{-1/p}\max_{j\leq 2t/\mu}\max_{k\leq \varepsilon t^{1/p}}|S(T_{j+k})-S(T_j)|= $$ $$=\limsup_{t\to\infty}t^{-1/p}\max_{j\leq 2t/\mu}\max_{k\leq \varepsilon t^{1/p}}|S(T_{j+k})-S(T_j)|=$$ $$ = \limsup_{t\to\infty}t^{-1/p}\max_{j\leq 2t/\mu}\max_{k\leq \varepsilon t^{1/p}}\Bigl|\varkappa k \mu  + v(B_{j+k}-B_j)+\beta\sqrt{\di \tau_1}(\widetilde{B}_{j+k}-\widetilde{B}_{j}) \Bigr|. $$
Due to the Cs\"{o}rg\H{o}-R\'{e}v\'{e}sz law for Wiener process \cite{CR} the last upper limit equals $\varepsilon \varkappa \mu, $ and since $\varepsilon$ was arbitrary, the desired relation holds. For $q=3,4,5,6$  the statement follows from \eqref{2}, \eqref{5},  \eqref{6} and \eqref{7a} respectively. For $q=7$ one should e.g. again use the Cs\"{o}rg\H{o}-R\'{e}v\'{e}sz law  noting that $|N(y)-u/\gamma|\leq 1,$ and finally, $\Phi_8$ is uniformly bounded. Theorem \ref{TheoremMain2}
is proved.

Let us now prove Theorem \ref{TheoremMain}.
From this point till Lemma \ref{RandomNagaev}, we always consider pairs $(t,x)$ such that 
\begin{equation}\label{pair}x \leq \frac{ t}{\log t}.\end{equation}
\begin{lm1}\label{Poisson}
We have $$\p\left(N^{-1}(t/\gamma)\geq 2t/\mu\right) \leq 2(e/2)^{  -x/\gamma },$$
for $(t,x)$ satisfying $\eqref{pair}.$
\end{lm1}

\begin{proof} $N^{-1}(t/\gamma)$ is the sum of $[t/\gamma]+1$ i.i.d. exponential random variables with parameter $\lambda,$ hence
$$\p(N^{-1}(  t/\gamma)\geq 2t/\mu)\leq  e^{-t/\gamma} 2^{[t/\gamma]+1},$$
 and it remains to use the inequality \eqref{pair}.
\end{proof}
   
 \begin{lm1}\label{MerlevedeRio}
There is  $a_5>0$ such that for $q=3,\ldots,8$ we have $$\p\left(\sup_{u\leq t}   |\Phi_q(u)|I\{N^{-1}(t/\gamma)\leq 2t/\mu \} > x \right) \leq a_5t x^{-p},$$
 for all $(t,x)$ satisfying $\eqref{pair}$.
\end{lm1}

\begin{proof}
For $q=3$ one should apply \eqref{3} and the estimate
\begin{equation}\label{Levy}\p(\sup_{u\leq t}|B_u-B_{[u]}|\geq x)\leq 4(t+1)e^{-x^2/2} .\end{equation}
For $q=4,5,6$ the corresponding statements directly follow from \eqref{4a}, \eqref{deheuvels}, \eqref{7} respectively.
For $ q=7$, using that $|u/\gamma - N(y)|\leq 1,$ the estimate follows analogously to \eqref{Levy}, and $\Phi_8$ is uniformly bounded.
\end{proof}

\begin{lm1}\label{Renewal}
There is  $a_6>0$ such that   $$\p\left(m(t) >\frac{2t}{\mu}\right) \leq a_6tx^{-p}, \,\, for\, all\,\, (t,x)\,\, satisfying\,\, \eqref{pair}.$$
\end{lm1}

\begin{proof}
Picking $b>0$ such that $\e\exp\{b(\mu/2-\tau_1)\}<1$ we have
$$\p \left(m(t) >\frac{2t}{\mu}\right) \leq
\p\left(- T_{[2t/\mu]}<=-t\right) \leq
e^{bt}(\e e^{-b\tau_1})^{[2t/\mu]} \leq 
\frac{1}{\e e^{-b\tau_1}} \left(  \e\exp\{b(\mu/2-\tau_1)\}\right)^{2x/\mu}.
$$

\end{proof}

 \begin{lm1}\label{Phi1}
There is  $a_7>0$ such that   $$\p\left(\sup_{u\leq t}   |\Phi_1(u)| \geq x \right) \leq a_7tx^{-p},\,\, for\, all\,\, (t,x)\,\, satisfying\,\, \eqref{pair}.$$
\end{lm1}

\begin{proof}
Note that $\sup_{u\leq t}|\Phi_1(u)|\leq \max_{k\leq m(t)+1}\eta_k,$ since $m(u)$ is the last renewal point of $T$ happening before $u,\,u\leq t.$ Consequently, by Lemma \ref{Renewal} we have 
$$\p\left(\sup_{u\leq t}|\Phi_1(u)|\geq x\right)\leq \p\left(m(t) >\frac{2t}{\mu}\right) +\p\left(\max_{k\leq 2t/\mu+1}\eta_k \geq    x\right)\leq    a_6tx^{-p} +
\frac{2t+\mu}{\mu}\e \eta_1^px^{-p} .$$
\end{proof}

To handle $\Phi_2$   we need the following two lemmas.

\begin{lm1}\label{OttavianiNagaev}
Suppose that $c>0$ and $X_1,X_2,\ldots$ are i.i.d. centered random variables such that $\e |X_1|^p<\infty$ for some $p>2,$ and let $Q_n=\sum_{j=1}^n X_j$ $(Q_0=0).$ Then  there exist $a_8=a_8(X_1,c)>0$  such that for any $n\in\mathbb{N}$ and all $x\geq c n^{1/p},\,x\leq n,$ one has
$$ \p\left(\max_{j\leq n}\max_{k\leq   x,j+k\leq n} (Q_{j+k} - Q_j)\geq  x\right) \leq  a_8nx^{-p}.$$
\end{lm1}

\begin{proof} We can assume that $[x]\geq 1. $ Divide the set $\{1,\ldots,n\}$ into consecutive blocks
$D_1,\ldots,D_M$ where all the blocks except for the last one have length $[x]$ and the last one is not longer than $[x].$ For $m\in\{1,\ldots,n\}$ let  $\varphi(m)=\max\{q[x]: q\in\mathbb{Z}_+,q[x]<m\}.$ 
We have
$$\p\left(\max_{j\leq n}\max_{k\leq   x, j+k\leq n} (Q_{j+k} - Q_j)\geq x\right)
=$$ $$=\p\left(\max_{j\leq n}\max_{k\leq   x, j+k\leq n} (Q_{j+k}-Q_{\varphi(j+k)}+
Q_{\varphi(j+k)}-Q_{\varphi(j)}+Q_{\varphi(j)}- Q_j)\geq  x\right)
\leq $$ $$\leq
\p\left(\max_{j=1,\ldots,M}\max_{k\leq x}|Q_k-Q_{\varphi(k)}|\geq\frac{x}{3}\right)\leq
M\p\left(\max_{k\leq x}|Q_k|\geq \frac{x}{3} \right)\leq 3\left(\frac{n}{[x]}+1\right)
\max_{k\leq x}\p\left(|Q_k|\geq \frac{x}{9}\right)$$
by the L\'{e}vy-Ottaviani inequality \cite[Prop. 1.1.2]{PenaGine}, and it remains to apply the Nagaev inequality \cite{Nagaev}.
\end{proof}

 \begin{lm1}\label{Phi2}
There is  $a_9 >0$ such that   $$\p\left(\sup_{u\leq t}  |\Phi_2(u)| \geq  x \right) \leq a_9tx^{-p},\,\, for\, all\,\, (t,x)\,\, satisfying\,\, \eqref{pair}.$$
  
\end{lm1}

\begin{proof}
$$\p\left(\sup_{u\leq t} |\Phi_2(u)| \geq  x \right)=
\p\left(\sup_{u\leq t}|S(T_{m(u)}) - S(T_{[y]})|\geq x\right)
\leq 
\p(N^{-1}(u/\gamma)>2t/\mu)+\p(m(t)>2t/\mu) +$$ $$+\p\left(\sup_{u\leq t}|m(u)-y(u)|>x-1\right)+\p\left(\sup_{k,m\leq 2t/\mu,|k-m|\leq x}|S(T_{k})-S(T_m)|\geq x\right). $$
First three summands are estimated by Lemma \ref{Poisson}, Lemma \ref{Renewal}
 and the relation \eqref{m-y} respectively. The fourth summand is estimated by Lemma \ref{OttavianiNagaev} for all $(t,x)$ such that $x\leq [2t/\mu] ,$ which is satisfied for all $t$ large enough due to \eqref{pair}.
\end{proof}

Consider now the case when $x $ lies in a domain of large deviations. We need the following analog of Nagaev inequality for random sums.

\begin{lm1}\label{RandomNagaev}
Let $X_1,X_2,\ldots$ be a sequence of centered random vectors such that $\{(X_n,\tau_n),n\geq 1\}$ are i.i.d. and $\e|X_1|^p<\infty,$ and let $Q_n=\sum_{j=1}^n X_j$ $(Q_0=0).$ Then there exists $a_{10}=a_{10}(X_1,\tau_1)>0$ depending on the distribution of $(X_1,\tau_1)$ such that for any $t\geq e$ and $x\geq t/\log t $ one has
$$\p\left(\max_{k\leq m(t)+1}|Q_k|>x\right)\leq a_{10} tx^{-p}. $$
\end{lm1}

\begin{proof}
Let $M_0\geq 1$ be a fixed integer such that $(\e e^{-\tau_1})^{M_0/2}<1/e.$ We have 
$$\p\left(\max_{k\leq m(t)+1}\left|Q_k\right|>x\right)\leq \p\left(\max_{k\leq M_0t+1}\left|Q_k\right|>x\right)+$$ $$+\sum_{M=M_0}^{\infty}x^{-p}\e \max_{k\leq (M+2)t}\left|Q_k\right|^pI\{Mt<m(t) \}=:R_0+x^{-p}\sum_{M=M_0}^{\infty}R_M. $$
Due to the L\'{e}vy-Ottaviani and Nagaev inequalities there exists $a_{11}>0$ such that $R_0\leq a_{11}tx^{-p}.$
To estimate $R_M$ ($M \geq M_0$) write, using H\"{o}lder inequality,
$$R_M\leq  \e \Biggl(\sum_{k\leq  (M+2)t}\left|X_k\right|\Biggr)^pI\{Mt<m(t)\} 
\leq ((M+2)t)^{p-1} \sum_{k\leq  (M+2)t}\e  |X_k|^pI\left\{- T_{[Mt]}\geq -t\right\}\leq
$$ $$ \leq ((M+2)t)^{p-1} \sum_{k\leq  (M+2)t}\e|X_k|^p \exp\left\{t - \sum_{i=1}^{[Mt]}\tau_i  \right\}.
$$
Note that for each $k$, the random vector $X_k$ is independent of all $\tau_i$'s ($i=1,\ldots,[Mt]$) except for at most one. Consequently,
$$\sum_{M=M_0}^{\infty}R_M\leq \sum_{M=M_0}^{\infty}(M+2)^{p}t^{p}  \e|X_1|^pe^t (\e \exp\{-\tau_1\})^{[Mt]} $$
which is a bounded function in $t\geq 1.$
\end{proof}

\begin{lm1}\label{LD} The relation $\eqref{1}$ holds for all $x > t/\log t.$ 
\end{lm1}

\begin{proof} First, by standard properties of Brownian motion one has
$$\p\left(\sup_{u\leq t}|W_u| > \frac{x}{2}\right)\leq 4de^{-x/(16\log x)}\leq a_{12}x^{-p} $$
for any $x$ and some absolute constant $a_{12}>0,$ provided that the pair $(x,t)$ satisfies Lemma's condition.
 Next, write
$$\p\left(\sup_{u\leq t}|S(u) - \varkappa u| > \frac{x}{2}\right)\leq \p\left(\sup_{k\leq m(t)}|S(T_k) - \varkappa T_k| > \frac{x}{6}\right) + $$ $$+\p\left(\sup_{u\leq t}|S(T_{m(u)}) - S(u)| > \frac{x}{6}\right) + \p\left(\sup_{u\leq t}|\varkappa||T_{m(u)} - u| > \frac{x}{6}\right) \leq $$ $$ \leq
\p\left(\sup_{k\leq m(t)}|S(T_k) - \varkappa T_k| > \frac{x}{6}\right) + \p\left(\sup_{k\leq m(t)+1}\eta_{k} > \frac{x}{6}\right) + \p\left(\sup_{k\leq m(t)+1}\tau_{k} > \frac{x}{6|\varkappa|}\right)  =:\sum_{q=1}^3 R_3.
$$
Lemma \ref{RandomNagaev} implies that $R_1\leq 6^p a_{10}(\xi_1-\varkappa \tau_1,\tau_1)tx^{-p}.$ Furthermore, let $t$ be large enough so that $x/12>\e\eta_1$ for  any $x\geq t/\log t.$
Then applying  Lemma \ref{RandomNagaev} to the sequence $\{\eta_k-\e\eta_k,\,k\geq 1\}$   one has
$$R_2\leq 12^pa_{10}(\eta_1-\e\eta_1,\tau_1) tx^{-p} ,  $$
and $R_3$ is estimated analogously.
\end{proof}

Theorem now follows from Lemmas \ref{Poisson}, \ref{MerlevedeRio}, \ref{Phi1}, \ref{Phi2}, \ref{LD}.

\end{proof}


\begin{thebibliography}{40}
  \bibitem{ASAExp} E.Bashtova, A.Shashkin. Strong Gaussian approximation for cumulative processes, arXiv:2006.09583.
 \bibitem{Berger} E.Berger. Fast sichere Approximationen von Partialsummen unabh\"{a}ngiger und
station\"{a}rer ergodischer Folgen von Zufallsvectoren, Dissertation, Universit\"{a}t zu
G\"{o}ttingen, 1982.
 \bibitem{CHS} M.Cs\"{o}rg\H{o}, L.Horv\'{a}th, J.Steinebach. Invariance principles for renewal processes. Ann. Probab., 1987, V. 15, N. 4, p. 1441--1460.
 \bibitem{CR} M.Cs\"{o}rg\H{o}, P.R\'{e}v\'{e}sz. How big are the increments of a Wiener process? Ann. Probab., 1979, V. 7, N. 4, p. 731--737. 
\bibitem{Einmahl} U.Einmahl. Extensions of results of Koml\'{o}s,  Major, and Tusn\'{a}dy to the multivariate case. J. Multival. Anal., 1989, V. 28, N. 1, p. 20--68.
\bibitem{EMK} P.Embrechts,  C.Kl\"{u}ppelberg, T.Mikosch. Modelling extremal events: for insurance and finance. Springer, 2013.
\bibitem{PenaGine} E.Gin\'{e}, V.de la Pe\~{n}a. Decoupling: from dependence to independence. Springer, 1999.
\bibitem{KMT} J.Koml\'{o}s,  P.Major, G.Tusn\'{a}dy. An Approximation of partial sums
of independent RV's, and the sample DF. I. Z. Wahrsch. verw. Geb., 1975, V. 32, N. 1, p. 111--131.
\bibitem{KMT2} J.Koml\'{o}s,   P.Major, G.Tusn\'{a}dy. An Approximation of partial sums
of independent RV's, and the sample DF. II. Z. Wahrsch. verw. Geb., 1976, V. 34, N. 1, p. 33--58.
\bibitem{Maj} P.Major. The approximation of partial sums of independent RV's. 
Z. Wahrsch. Verw. Geb., 1976, V. 35, p. 213--220.
 \bibitem{MerRio} F.Merlev\'{e}de, E.Rio. Strong approximation for additive functionals of geometrically ergodic Markov chains. Electron. J. Probab., 2015, V. 20, N. 14, p. 1--27.
 \bibitem{Nagaev} S.V.Nagaev. Some limit theorems for large deviations. Theory Probab. Appl., 1965, V. 10, N. 2, p. 214--235.
\bibitem{Smith1955} W.L.Smith. Regenerative stochastic processes. Proc. Royal Soc. London Ser. A, 1955, V. 232, N. 1188, p. 6--31. 
 \bibitem{Strass} V.Strassen. An invariance principle for the law of the iterated logarithm. Z. Wahrsch. Verw. Geb., 1964, V. 3,  p. 211--226.
\bibitem{Zait_survey} A.Yu.Zaitsev. The accuracy of strong Gaussian approximation for sums of independent random vectors. Russ. Math.Surv., 2013, V. 68, N.  4, pp.721--761.

 
 \end{thebibliography}
\end{document}